\documentclass{amsart}

\usepackage{enumerate}

\newcommand{\R}{\mathbb{R}}

\newcommand{\N}{\mathbb{N}}
\newcommand{\hd}{\dim_{\textup{H}}}

\newcommand{\ubd}{\overline{\dim}_{\textup{B}}}

\newcommand{\Lip}{\textnormal{Lip}}
\newcommand{\si}{S_{\mathbf{i}}}

\newcommand{\I}{\mathcal{I}}

\newcommand{\hm}{\mathcal{H}}

\newtheorem{theorem}{Theorem}[section]
\newtheorem{lemma}[theorem]{Lemma}
\newtheorem{cor}[theorem]{Corollary}

\theoremstyle{definition}

\theoremstyle{remark}

\numberwithin{equation}{section}
\allowdisplaybreaks

\begin{document}
\title[On the dimension and measure of inhomogeneous attractors]{On the dimension and measure of inhomogeneous attractors}
\author{Stuart A. Burrell}
\address{S.A. Burrell, School of Mathematics and Statistics, University of St Andrews, St Andrews, KY16 9SS, United Kingdom.}
\email{sb235@st-andrews.ac.uk}


\date{\today}
\dedicatory{}

\begin{abstract}
A central question in the field of inhomogeneous attractors has been to relate the dimension of an inhomogeneous attractor to the condensation set and associated homogeneous attractor. This has been achieved only in specific settings, with notable results by Olsen, Snigireva, Fraser and K\"aenm\"aki on inhomogeneous self-similar sets, and by Burrell and Fraser on inhomogeneous self-affine sets. This paper is devoted to filling a significant gap in the dimension theory of inhomogeneous attractors, by studying those formed from arbitrary bi-Lipschitz contractions. We show that the maximum of the dimension of the condensation set and a quantity related to pressure, which we term upper Lipschitz dimension, forms a natural and general upper bound on the dimension. Additionally, we begin a new line of enquiry; the methods developed are used to investigate the Hausdorff measure of inhomogeneous attractors. Our results have applications for affine systems with affinity dimension less than or equal to one and systems satisfying bounded distortion, such as conformal systems in dimensions greater than one. 
\end{abstract}
\maketitle

\section{Introduction}
Let $(X, d)$ be a compact metric space. A map $S : X \rightarrow X$ is a contraction on $X$ if there exists  $c \in (0, 1)$ such that 
\begin{equation*}
d(S(x), S(y)) \leq c d(x,y)
\end{equation*}
for all $x, y \in X$, and a similarity with ratio $c$ if $\leq$ can be replaced with $=$. We call a finite collection $\mathbb{I} = \{S_i\}_{i = 1}^{N}$ of contractions on $X$ an \emph{iterated function system} (IFS). It is well known that for each IFS there exists a unique non-empty compact set $F$, called the \emph{homogeneous} attractor, such that
\begin{equation*}
F = \bigcup\limits_{i = 1}^{N}S_i(F).
\end{equation*}
Analogously, if we fix a compact set $C \subseteq X$, there exists a unique non-empty compact set $F_C$ such that
\begin{equation*}
F_C = \bigcup\limits_{i = 1}^{N}S_i(F_C) \cup C,
\end{equation*}
called the \emph{inhomogeneous} attractor with \emph{condensation set} $C$. This construction was introduced by Barnsley (1985) \cite{barn} and Hata (1985) \cite{hata}, and has been studied in, for example, \cite{fraserbaker,burrellfraser,Fraser,18,antti,olssni,sniphd}. For an illustration, see \cite{Fraser}.

To proceed further we require some notation. For $\mathbb{I} = \{S_i\}_{i = 1}^{N}$, define $\mathcal{I} = \{1,\dots,N\}$ and write $\si = S_{i_1} \circ \dots \circ S_{i_k}$ where $\textbf{i} = (i_1,i_2,\dots,i_k) \in \mathcal{I}^k$. Furthermore, let
\begin{equation*}
\mathcal{I}^* = \bigcup\limits_{k = 1}^{\infty} \mathcal{I}^k
\end{equation*}
denote the set of finite words over $\mathcal{I}$.
An elegant and useful expression for $F_C$, from \cite{Fraser} and \cite{sniphd}, is
\begin{equation*}
F_C = F_{\emptyset} \cup \mathcal{O},
\end{equation*}
where $F_{\emptyset}$ is the corresponding homogeneous attractor of $\mathbb{I}$ and $\mathcal{O}$ is the \emph{orbital} set defined by
\begin{equation*}
\mathcal{O} = C \cup \bigcup\limits_{\textbf{i} \in \mathcal{I}^*} \si(C).
\end{equation*}
For a given dimension, it is natural to ask whether, or in what situations, 
\begin{equation}\label{target1}
\dim F_C = \max\{\dim F_\emptyset, \dim C\},
\end{equation}
where $\dim$ denotes some notion of dimension. If the dimension is known to be countably stable, this relationship is immediate. In \cite{Fraser}, an answer is given for the countably unstable upper box dimension in the case where $\mathbb{I}$ consists only of similarity mappings. Recall that for a non-empty set $F\subseteq X$, the upper box dimension is defined as
\begin{equation*}
\overline{\dim}_B F = \limsup_{\delta \rightarrow 0} \frac{\log N_{\delta}(F)}{-\log \delta },
\end{equation*}
where $N_{\delta}(F)$ denotes the minimum number of balls of diameter $\delta$ required to cover $F$. In the case of the analogously defined lower box dimension, Fraser shows (\ref{target1}) does not hold in general for similarities \cite{Fraser}. Moreover, even for upper box dimension the relationship may not hold if the IFS has significant overlaps \cite{fraserbaker}.\\

The works of Fraser \cite{Fraser}, K\"aenm\"aki \cite{antti}, Olsen and Snigireva \cite{olssni}, and Snigireva \cite{sniphd} deal exclusively with inhomogeneous self-similar sets, and rely on the notion of similarity dimension. Recall that, for a finite collection of contracting similarities $\{S_i\}_{i = 1}^{N}$, the similarity dimension is the value of $s$ satisfying Hutchinson's formula
\begin{equation*}
\sum\limits_{i = 1}^{N} \Lip(S_i)^s = 1,
\end{equation*}
where $\Lip(S_i)$ denotes the similarity ratio of $S_i$. Separation conditions are needed to guarantee equality of the upper box and similarity dimensions, see \cite{Falconer} for details. One such condition is the \emph{strong open set condition} (SOSC), which requires that there exists a non-empty open set $U$ such that $F\cap U \neq \emptyset$ and
\begin{equation*}
\bigcup\limits_{i = 1}^{N} \si(U) \subseteq U,
\end{equation*}
with the union disjoint. It is known from \cite[Theorem 2.1]{Fraser} that an inhomogeneous self-similar set $F_C$ satisfies 
\begin{equation*}
\max\{\ubd F_{\emptyset},\, \ubd C\} \leq \ubd F_{C}
\leq \max\{s,\ \ubd C\},
\end{equation*}
where $s$ denotes the similarity dimension. Thus, if the SOSC is satisfied, then we obtain equality above and recover the desired relationship (\ref{target1}). To obtain such a bound for general classes of maps we construct an analogue of similarity dimension. For $S: X \rightarrow X$, let
\begin{equation*}
\Lip^+(S) = \sup\limits_{\substack{x,y \in X\\x \neq y}} \frac{d(S(x), S(y))}{d(x, y)}
\end{equation*}
and
\begin{equation*}
\Lip^-(S) = \inf\limits_{\substack{x,y \in X\\x \neq y}} \frac{d(S(x), S(y))}{d(x, y)}
\end{equation*}
denote the upper and lower Lipschitz constants respectively. We say $S$ is a bi-Lipschitz contraction if $0 < \Lip^-(S) \leq \Lip^+(S) < 1$. Next, for an IFS $\mathbb{I} = \{S_i\}_{i = 1}^{N}$, define $(s_k)$ to be the solution of
\begin{equation}\label{skdef1}
\sum\limits_{\textbf{i} \in \mathcal{I}_k} \Lip^+(\si)^{s_k} = 1
\end{equation}
for each $k \in \N$. The corresponding limit, given by
\begin{equation}\label{skdef2}
s = \lim_{k \rightarrow \infty} s_k,
\end{equation}
is referred to as the \emph{upper Lipschitz dimension}. The existence of the upper Lipschitz dimension is established by considering the pressure function $P(t) = \lim\limits_{k \rightarrow \infty} P_k(t)$ where
\begin{equation*}
P_k(t) = \frac{1}{k}\log \sum\limits_{\textbf{i} \in \mathcal{I}^k} \Lip^+(\si)^t.
\end{equation*}
Subadditivity and Fekete's lemma imply $P(t)$ exists for all $t \geq 0$. Moreover, it is well known $P$ is continuous, monotonically decreasing and has a unique zero. Since $P_k \rightarrow P$ pointwise, it follows that the upper Lipschitz dimension exists and is equal to the zero of $P$. For further details on pressure we direct the reader to \cite[Chapter 5]{Techniques} and the references therein.\\

Our main result establishes bounds on the upper box dimension of $F_C$ for general IFSs consisting of bi-Lipschitz contractions. The methodology of Fraser relies heavily on the multiplicativity of the pressure function for similarities, which presents complications in the general case. To overcome this, we show that dimension is invariant under passing to a derived system with desirable properties, see Lemma \ref{compattr} and Corollary \ref{keycore}.\\

Additionally, we begin a new line of enquiry; the methods developed are used to investigate the Hausdorff measure of inhomogeneous attractors. This may serve as a platform for future work, such as deeper investigations into the sensitive dependence of the measure with separation conditions on $C$. \\

We also present applications where further assumptions yield stronger corollaries in popular contexts, such as conformal and low dimensional affine systems. Hereafter, definitions relevant only to an individual application or result are contained in the corresponding section.

\section{Results}
\subsection{Dimensions of inhomogeneous attractors}
The following provides an analogue of Fraser's result on inhomogeneous self-similar sets \cite{Fraser} for IFSs consisting of arbitrary bi-Lipschitz maps.
\begin{theorem}\label{main2}
Let $(X,d)$ be a compact metric space and $\mathbb{I} = \{S_i\}_{i = 1}^{N}$ denote an IFS consisting of bi-Lipschitz maps with compact condensation set $C \subseteq X$. We have
\begin{equation*}
\max\{\overline{\dim}_B F_{\emptyset},\, \overline{\dim}_B C\} \leq \overline{\dim}_B F_C \leq \max\left\{{s, \overline{\dim}_B C}\right\},
\end{equation*}
where $s$ is equal to the upper Lipschitz dimension. 
\end{theorem}

Theorem \ref{main2} implies that if $\ubd C \geq s$, then $\ubd F_C = \ubd C$. It is therefore natural to consider applications in which $s = \ubd F_\emptyset$. One such scenario involves the notion of \emph{bounded distortion}. An IFS $\mathbb{I} = \{S_i\}_{i = 1}^{N}$ satisfies the property of bounded distortion if there exists some uniform constant $L > 1$ such that
\begin{equation*}
\frac{\Lip^+(\si)}{\Lip^-(\si)} < L,
\end{equation*}
for all $\textbf{i} \in \mathcal{I}^*$. Lemma \ref{supperbound} and a simple modification of \cite[Proposition 9.7]{Falconer} imply that bounded distortion together with the SOSC force $s = \ubd F_\emptyset$. This immediately yields the following corollary of Theorem \ref{main2}.

\begin{cor}\label{bd}
Let $(X,d)$ be a compact metric space and $\mathbb{I} = \{S_i\}_{i = 1}^{N}$ denote an IFS satisfying bounded distortion with compact condensation set $C \subseteq X$. If $\mathbb{I}$ satisfies the SOSC, then
\begin{equation*}
\overline{\dim}_B F_C = \max\left\{{\ubd{F_\emptyset}, \overline{\dim}_B C}\right\}.
\end{equation*}
\end{cor}

In $\R^n$ for $n \geq 2$, the conformal maps (that is, locally angle-preserving maps) are a well-known example that satisfy bounded distortion, as well as $C^{1+\alpha}$ maps on $\R$ \cite{Techniques}. Thus, conformal iterated function systems, formally defined and studied in such papers as \cite{Feng}, constitute a useful application of this result.\\

Finally, Theorem \ref{main2} provides an extremely succinct proof of the main result of \cite{burrellfraser} for low-dimensional affine systems. For an introduction to the theory of self-affine sets, see \cite{burrellfraser} or \cite{Falconer}. If the affinity dimension is less than or equal to one, then it coincides with the upper Lipschitz dimension, since $\Lip^+(S)$ corresponds to the largest singular value of the linear component of $S$. This yields the following.

\begin{cor}\label{bd1}
Let $\mathbb{I} = \{S_i\}_{i = 1}^{N}$ be an affine IFS with compact condensation set $C \subseteq X$ and affinity dimension $s$. If $s \leq 1$, then
\begin{equation*}
\max\left\{{\ubd{F_\emptyset}, \overline{\dim}_B C}\right\} \leq \overline{\dim}_B F_C \leq \max\left\{{s, \overline{\dim}_B C}\right\}.
\end{equation*}
\end{cor}
In particular, Corollary \ref{bd1} implies that if the affinity dimension is less than or equal to one and equals $\ubd F_\emptyset$, then (\ref{target1}) is satisfied. Falconer shows in \cite{1988kjf} that the affinity and upper box dimensions coincide almost surely upon randomizing the translations, even if the SOSC fails, and it follows from recent breakthrough results of B\'ar\'any, Hochman and Rapaport \cite{bhr} that mild assumptions are sufficient to force equality in the plane. For a more detailed discussion, see \cite{burrellfraser}. However, it is worth noting that (\ref{target1}) does not always hold in the affine setting. In particular, from the results of Fraser \cite{18}, it is possible to construct simple examples of inhomogeneous self-affine sets with affinity dimension $s < 1$ satisfying
\begin{equation*}
\max\left\{{\ubd{F_\emptyset}, \ubd C}\right\} < 
\ubd F_C < \max\left\{{s, \ubd C}\right\}.
\end{equation*}

\subsection{Hausdorff measure of inhomogeneous attractors}
The approach developed to study dimension, specifically the strategy implied by Corollary \ref{keycore}, allows us to investigate the Hausdorff measure of $F_C$ at the critical value. Recall that, for $0 \leq t \leq n$ and $\delta > 0$, the $t$-dimensional $\delta$-approximate Hausdorff measure of $F \subset \R^n$ is 
\begin{equation*}
\mathcal{H}^t_\delta(F) = \inf\left\{\sum\limits_{i = 1}^{\infty}|U_i|^t :F\subset \bigcup\limits_{i = 1}^{\infty} U_i, 0 < |U_i| < \delta \right\},
\end{equation*}
and 
\begin{equation*}
\mathcal{H}^t(F) = \lim\limits_{\delta \rightarrow 0} \mathcal{H}^t_\delta(F)
\end{equation*}
is the $t$-dimensional Hausdorff measure. This leads to the definition of Hausdorff dimension as
\begin{equation*}
\dim_H F = \inf\{t: \mathcal{H}^t(F) = 0\} = \sup\{t:\mathcal{H}^t(F) = \infty\},
\end{equation*}
with $t = \dim_H F$ known as the critical value. It has been of historical interest (e.g. \cite{Falconer}) to compute $\mathcal{H}^t(F)$ at the critical value. It may well be zero, finite or infinite. Recall that countable stability and monotonicity of Hausdorff dimension \cite{Falconer} readily imply that
\begin{equation*}
\hd F_C = \max\{ \hd F_{\emptyset}, \hd C \},
\end{equation*}
and we investigate $\mathcal{H}^t(F_C)$ in each case. Although we omit the details, the following results readily extend to any family of measures satisfying the scaling property and their associated dimension, such as packing measures.\\

For $t = \hd F_C$, if $\hm^t(C) = 0$ or $\hm^t(C) = \infty$, then it is clear that $\hm^t(F_C) = \hm^t(F_\emptyset)$ and $\hm^t(F_C) = \infty$, respectively. Thus, our main theorem deals with the case where $\hm^t(C)$ is positive and finite. For the problem to be tractable, some separation conditions are required. A natural choice is the \emph{condensation open set condensation}, a modification of the SOSC adapted for the inhomogeneous case, as utilised in \cite{burrellfraser,antti,olssni,sniphd}. An IFS satisfies the COSC if there exists an open set $U$ with 
$$C \subset U \setminus \bigcup\limits_{i = 1}^{N} \overline{S_i(U)},$$ such that $S_i(U)\subset U$ for $i = 1,\dots,N$, and
$i \neq j \implies S_i(U) \cap S_j(U) = \emptyset$.

\begin{theorem}\label{hausthm}
Let $(X,d)$ be a compact metric space and $\mathbb{I} = \{S_i\}_{i = 1}^{N}$ denote an IFS with compact condensation set $C \subseteq X$ and upper Lipschitz dimension $s$. Suppose $t = \ubd F_C$ and $0 < \mathcal{H}^t(C) < \infty$. It follows that
\begin{enumerate}[i)]
\item if $t > s$, then $0 < \mathcal{H}^t(F_C) < \infty$;
\item if $\mathbb{I}$ satisfies the COSC, then
$$
\hm^t(F_C) \geq \mathcal{H}^t(F_\emptyset) + \mathcal{H}^t(C)\left(1 + \sum\limits_{k = 1}^{\infty}\left(\sum\limits_{i \in \mathcal{I}} \Lip^-(S_i)^t\right)^k\right)\hphantom{;}
$$
and
$$
\hm^t(F_C) \leq \mathcal{H}^t(F_\emptyset) + \mathcal{H}^t(C)\left(1 + \sum\limits_{k = 1}^{\infty}\left(\sum\limits_{i \in \mathcal{I}} \Lip^+(S_i)^t\right)^k\right).
$$
\end{enumerate}
\end{theorem}

This yields a pleasing closed form expression for inhomogeneous self-similar sets, as studied in, for example, \cite{Fraser,antti,sniphd}.
\begin{cor}\label{hauscor}
Let $(X,d)$ be a compact metric space and $\mathbb{I} = \{S_i\}_{i = 1}^{N}$ denote an IFS consisting of similarities satisfying the COSC with compact condensation set $C \subseteq X$ and similarity dimension $s$. If $t = \hd F_C > s$ and $0 < \hm^t(C) < \infty$, then 
\begin{equation*}\label{closed}
\mathcal{H}^t(F_C) = \frac{\mathcal{H}^t(C)}{1 - \sum\limits_{i \in \mathcal{I}} \Lip(S_i)^t}.
\end{equation*}
\end{cor}

\begin{proof}
For a similarity $S$, we have $\Lip^+(S) = \Lip^-(S)$, and the result follows immediately from Theorem \ref{hausthm} (ii), since the upper Lipschitz and similarity dimensions coincide.
\end{proof}

We hope the above may serve as a platform for future developments. In particular, it would be interesting to discover alternative conditions to the COSC that control the sensitive interaction between $F_\emptyset$ and $\mathcal{O}$ while yielding similar results.

\section{Proof of Theorem \ref{main2}} \label{mainproof}

Let $\mathbb{I} = \{S_i\}_{i=1}^{N}$ be an IFS consisting of bi-Lipschitz maps and $C \subseteq X$ be compact.

\subsection{Preliminary lemmas}
\begin{lemma}\label{compattr}
For all $k \in \N$, the IFS given by $\mathbb{I}_k = \{\si\}_{\mathbf{i} \in \mathcal{I}^k}$ satisfies
\begin{equation*}
\ubd F_C = \ubd F_C^k,
\end{equation*}
where $F_C^k$ denotes the inhomogeneous attractor associated with $\mathbb{I}_k$.
\begin{proof}
Fix $k \in \N$ and observe that
\begin{equation*}
F_\emptyset = \bigcup\limits_{i \in \mathcal{I}} S_i(F_\emptyset) = \bigcup\limits_{i_1 \in \mathcal{I}}\bigcup\limits_{i_2 \in \mathcal{I}} S_{i_1}(S_{i_2}(F_\emptyset)) = \dots = \bigcup\limits_{\textbf{i} \in \mathcal{I}^k} \si(F_\emptyset)
\end{equation*}
and so $F_\emptyset = F_\emptyset^k$, where $F_\emptyset^k$ denotes the unique homogeneous attractor associated with $\mathbb{I}_k$. Recall that $F_C = F_{\emptyset} \cup \mathcal{O}$ and $\mathcal{O} = C \cup \bigcup\limits_{\textbf{i} \in \mathcal{I}^*} \si(C)$. Hence,
\begin{equation*}
F_C = F_{\emptyset} \cup C \cup \bigcup\limits_{\textbf{i} \in \mathcal{I}^*} \si(C) 
\end{equation*}
and
\begin{equation*}
F_C^k = F_{\emptyset} \cup C\cup\bigcup\limits_{\textbf{i} \in (\mathcal{I}^k)^*} \si(C),
\end{equation*}
where $(\mathcal{I}^k)^*$ denotes all finite concatenations of length $k$ words over $\mathcal{I}$. Thus, by finite stability of box dimension, it suffices to show
\begin{equation}\label{target}
\ubd  \bigcup\limits_{\textbf{i} \in \mathcal{I}^*} \si(C) \leq \ubd  \bigcup\limits_{\textbf{i} \in (\mathcal{I}^k)^*} \si(C),
\end{equation}
since the opposite inequality follows immediately by monotonicity. Observe
\begin{equation*}\label{keyeq}
\ubd  \bigcup\limits_{\textbf{i} \in \mathcal{I}^*} \si(C) = \max\limits_{t = 1, \dots, k}\ubd  \bigcup\limits_{\substack{\textbf{i} \in \mathcal{I}^*\\ |\textbf{i}| = nk + t \\ n \geq 0}} \si(C),
\end{equation*}
and let $m$ be the value of $t$ that realises the maximum. First, note that
\begin{equation}\label{cbound}
\ubd C \leq \ubd \si(C) \leq \ubd \bigcup\limits_{\textbf{i} \in (\mathcal{I}^k)^*} \si(C)
\end{equation}
for all $\textbf{i} \in \mathcal{I}^*$, since $\si$ is bi-Lipschitz. Hence
\begin{align*}
\ubd  \bigcup\limits_{\textbf{i} \in \mathcal{I}^*} \si(C) &=\ubd  \bigcup\limits_{n = 0}^{\infty}\bigcup\limits_{\substack{\textbf{i} \in \mathcal{I}^*\\ |\textbf{i}| = nk + m}} \si(C) \\
&= \ubd \bigcup\limits_{\textbf{j} \in \mathcal{I}^m}\left(S_{\textbf{j}}(C)\cup\bigcup\limits_{\textbf{i} \in (\mathcal{I}^k)^*} (S_{\textbf{ji}}(C))\right) \\
&= \ubd \bigcup\limits_{\textbf{j} \in \mathcal{I}^m}S_{\textbf{j}}\left(C \cup \bigcup\limits_{\textbf{i} \in (\mathcal{I}^k)^*} (S_{\textbf{i}}(C))\right) \\
&= \max\limits_{\textbf{j} \in \mathcal{I}^m} \ubd  S_{\textbf{j}}\left(C \cup\bigcup\limits_{\textbf{i} \in (\mathcal{I}^k)^*}(S_{\textbf{i}}(C))\right) \\
&\leq \max\{\ubd C, \ubd \bigcup\limits_{\textbf{i} \in (\mathcal{I}^k)^*} \si(C)\}\\
&\leq \ubd \bigcup\limits_{\textbf{i} \in (\mathcal{I}^k)^*} \si(C)
\end{align*}
by (\ref{cbound}).
\end{proof}
\end{lemma}

This has the following corollary that is fundamental to our approach.

\begin{cor}\label{keycore}
For $t > \max\{s, \ubd C\}$, there exists a $K \in \N$ such that $t > s_k$ for all $k > K$, and each IFS given by $\mathbb{I}_k = \{\si\}_{\mathbf{i} \in \mathcal{I}^k}$ satisfies
\begin{equation*}
\ubd F_C = \ubd F_C^k.
\end{equation*}
\begin{proof}
Since $s_k \rightarrow s$, there exists $K \in \N$ such that $|s - s_k | \leq \frac{t - s}{2}$ for all $k > K$. The result then follows immediately from Lemma \ref{compattr}.
\end{proof}
\end{cor}

The next Lemma is analogous to \cite[Lemma 3.2]{Fraser} and illustrates the motivation for Corollary $\ref{keycore}$.

\begin{lemma}\label{btdef}
If $t > s_1$, then there exists a constant $b_t$ such that
\begin{equation*}
\sum\limits_{\mathbf{i} \in \mathcal{I}^*} \Lip^+(\si)^t = b_t < \infty.
\end{equation*}
\begin{proof}
Observe that $t > s_1$ implies
\begin{equation*}
\sum\limits_{i \in \mathcal{I}}\Lip^+(S_i)^t < 1.
\end{equation*}
Hence
\begin{align*}
\sum\limits_{\textbf{i} \in \mathcal{I}^*}\Lip^+(\si)^t &= \sum_{k = 1}^{\infty}\sum\limits_{\textbf{i} \in \mathcal{I}^k}\Lip^+(\si)^t\\ &\leq \sum_{k = 1}^{\infty}\left(\sum\limits_{i \in \mathcal{I}}\Lip^+(S_i)^t\right)^k \\&< \infty,
\end{align*}
by convergence of the geometric series.
\end{proof}
\end{lemma}

A natural way to construct efficient $\delta$-covers is to consider the finite set of cylinders $\si(X)$ such that $\Lip^+(\si) < \delta$ and $\Lip^+(S_{\textbf{i}_p}) \geq \delta$ for any prefix $\textbf{i}_p$ of $\textbf{i}$. For $\textbf{i} = (i_1,...,i_k) \in \mathcal{I}^*$ we let $\textbf{i}_- = (i_1,...,i_{k-1})$ and write $|\textbf{i}|$ to denote the length of the string $\textbf{i}$. If $\delta \in (0, 1]$, define the $\delta$-stopping, denoted $\mathcal{I}(\delta)$, by
\begin{equation*}
\mathcal{I}(\delta) = \{\textbf{i} \in \mathcal{I}^* : \Lip^+(\si) < \delta \leq \Lip^+(S_{\textbf{i}_-}) \}.
\end{equation*}
We assume for convenience that $\Lip^+(S_\omega) = 1$, where $\omega$ denotes the empty word. If $\textbf{i} \in \mathcal{I}^*$ satisfies $\Lip^+(\si) < \delta$, then it is clear there exists a prefix $\textbf{i}_p$ such that $\textbf{i}_p \in \mathcal{I}(\delta)$. To establish a bound on $|\mathcal{I}(\delta)|$, we define 
\begin{equation*}
L_{\min} = \min\limits_{i \in \mathcal{I}} \Lip^-(S_i) > 0.
\end{equation*}
\begin{lemma}\label{dtconst}
If $t > s_1$, then
\begin{equation*}
|\mathcal{I}(\delta)| \leq b_t L_{\min}^{-t}\delta^{-t}
\end{equation*}
for all $\delta \in (0, 1]$.
\begin{proof}
For $\textbf{i} \in \mathcal{I}(\delta)$, we have
\begin{equation}
\Lip^+(\si) \geq \Lip^+(S_{\textbf{i}_-})L_{\min} \geq \delta L_{\min} > 0.
\end{equation}
Hence
\begin{align*}
b_t \geq \sum\limits_{\textbf{i} \in \mathcal{I}(\delta)} \Lip^+(\si)^t \geq \sum\limits_{\textbf{i} \in \mathcal{I}(\delta)} (\delta L_{\min})^t = |\mathcal{I}(\delta)| (\delta L_{\min})^t 
\end{align*}
and the desired inequality follows immediately.
\end{proof}
\end{lemma}

This yields an alternative and succinct proof of the well-known result that the dimension of the homogeneous attractor is bounded above by the upper Lipschitz dimension.

\begin{lemma}\label{supperbound}
$\ubd F_\emptyset \leq s$, where $s$ denotes the upper Lipschitz dimension of $\mathbb{I}$.
\begin{proof}
Fix $\delta \in (0, 1]$ and let $t > s$ be arbitrary. By Corollary \ref{keycore}, we may assume $t > s_1$ while preserving $\ubd F_\emptyset$. The result then follows from Lemma \ref{dtconst} since the cylinder sets $\{\si(X) : \textbf{i} \in \mathcal{I}(\delta)\}$ form a $\delta$-cover of $F_\emptyset$, and so $N_\delta(F_\emptyset) \leq |\mathcal{I}(\delta)| \leq b_tL_{\min}^{-t}\delta^{-t}$.
\end{proof}
\end{lemma}

For clarity in our later calculation, we provide one further lemma.

\begin{lemma}\label{inclusionlemma}
For $\delta \in (0, 1]$, we have
\begin{equation*}
\bigcup\limits_{\substack{\mathbf{i} \in \mathcal{I}^* \\ \Lip^+(S_{\mathbf{i}}) < \delta}} S_{\mathbf{i}}(C) \subseteq \bigcup\limits_{\mathbf{i} \in \mathcal{I}(\delta)} S_{\mathbf{i}}(X).
\end{equation*}
\begin{proof}
If 
\begin{equation*}
x \in \bigcup\limits_{\substack{\textbf{i}\in \mathcal{I}^* \\ \Lip^+(S_{\textbf{i}}) < \delta}} S_{\textbf{i}}(C)
\end{equation*} there exists some $\textbf{i} = (i_1, i_2,\dots,i_n)  \in \mathcal{I}^*$ with $\Lip^+(\si) < \delta$ and a $c \in C$ such that $x = \si(c)$. Let $\textbf{i}_p = (i_1,i_2,\dots,i_p)$ denote the prefix of $\textbf{i}$ with $\textbf{i}_p \in \mathcal{I}(\delta)$, then $x = S_{\textbf{i}_p}(S_{(i_{p+1}, i_{p+2},\dots, i_{n})}(c)) \in S_{\textbf{i}_p}(X)$, as required.
\end{proof}
\end{lemma}

\subsection{Proof of Theorem \ref{main2}}

Monotonicity of upper box dimension implies
\begin{equation*}
\max\{\overline{\dim}_B F_{\emptyset},\, \overline{\dim}_B C\} \leq  \overline{\dim}_B F_{C},
\end{equation*}
since $F_{\emptyset} \cup C \subseteq F_{\emptyset} \cup \mathcal{O} = F_{C}$. Moreover, by finite stability of upper box dimension we have 
\begin{equation*}
\overline{\dim}_B F_{C} \leq \max\{ \overline{\dim}_B F_{\emptyset},\, \overline{\dim}_B  \mathcal{O} \}.
\end{equation*}
Hence, since $\ubd F_\emptyset \leq s$ by Lemma \ref{supperbound}, it suffices to show
\begin{equation*}
\overline{\dim}_B  \mathcal{O} \leq \max \{ s, \overline{\dim}_B  C \},
\end{equation*}
where $s$ denotes the upper Lipschitz dimension.\\

Let $t > \max\{s, \ubd C\}$. By Corollary \ref{keycore}, since our interest is in $\ubd F_C$, we can assume hereafter that $t > s_1$ without loss of generality. The definition of box dimension implies that there exists a constant $c_t$ such that
\begin{equation}\label{cconst}
N_{\delta}(C) \leq c_t\delta^{-t},
\end{equation}
for all $\delta \in (0, 1]$. Further, since $X$ is compact, $N_1(X)$ is a finite constant that does not depend on $t$. Observe
\begin{flalign*}
N_{\delta}(\mathcal{O}) &= N_{\delta}\left(C \cup \bigcup\limits_{\textbf{i} \in \mathcal{I}^*}  S_{\textbf{i}}(C)\right) \\
& \leq N_{\delta}(C) + N_{\delta}\left(\bigcup\limits_{\substack{\textbf{i} \in \mathcal{I}^* \\ \Lip^+(S_{\textbf{i}}) \geq \delta}} S_{\textbf{i}}(C)\right) + 
N_{\delta}\left(\bigcup\limits_{\substack{\textbf{i} \in \mathcal{I}^* \\ \Lip^+(S_{\textbf{i}}) < \delta}} S_{\textbf{i}}(C)\right) \\
& \leq N_{\delta}(C) + \sum\limits_{\substack{\textbf{i} \in \mathcal{I}^* \\ \Lip^+(S_{\textbf{i}}) \geq \delta}} N_{\delta}(S_{\textbf{i}}(C)) + 
N_{\delta}\left(\bigcup\limits_{\textbf{i} \in \mathcal{I}(\delta)} S_{\textbf{i}}(X)\right) \textnormal{\,\,\,\,(by Lemma \ref{inclusionlemma})}\\
& \leq N_{\delta}(C) + \sum\limits_{\substack{\textbf{i} \in \mathcal{I}^* \\ \Lip^+(S_{\textbf{i}}) \geq \delta}} N_{\delta}(S_{\textbf{i}}(C)) + 
\sum\limits_{\textbf{i} \in \mathcal{I}(\delta)} N_{\delta}(S_{\textbf{i}}(X))\\
& \leq N_{\delta}(C) + \sum\limits_{\substack{\textbf{i} \in \mathcal{I}^* \\ \Lip^+(S_{\textbf{i}}) \geq \delta}} N_{{\delta}/{\Lip^+(S_{\textbf{i}})}}(C) + 
\sum\limits_{\textbf{i} \in \mathcal{I}(\delta)} N_{\delta/\Lip^+(S_{\textbf{i}})}(X) \\
& \leq c_t\delta^{-t} + \sum\limits_{\substack{\textbf{i} \in \mathcal{I}^* \\ \Lip^+(S_{\textbf{i}}) \geq \delta}} c_t(\delta /\Lip^+(S_{\textbf{i}}))^{-t} + 
\sum\limits_{\textbf{i} \in \mathcal{I}(\delta)} N_1(X) \\
& \leq c_t\delta^{-t} + c_t\delta^{-t}\sum\limits_{\substack{\textbf{i} \in \mathcal{I}^*}} \Lip^+(S_{\textbf{i}})^{t} + 
|\mathcal{I}(\delta)|N_1(X) \\
&\leq c_t\delta^{-t} + c_t\delta^{-t}b_t + 
b_t L_{\min}^{-t}\delta^{-t}N_1(X) \textnormal{\,\,\,\,(by Lemmas \ref{btdef} and \ref{dtconst})}\\
&\leq \delta^{-t}( c_t + c_tb_t + 
b_t L_{\min}^{-t}N_1(X)).
\end{flalign*}
\hfill $\square$

\section{Proof of Theorem \ref{hausthm}}
Let $\mathbb{I} = \{S_i\}_{i =1}^{N}$ be an IFS and $C \subseteq X$ be compact.\\
 
We require two technical Lemmas which provide a similar strategy in this context as Lemma \ref{keycore} allowed for dimension.
\subsection{Preliminary lemmas}

\begin{lemma}\label{hauskey}
Let $t \geq 0$ and suppose $\mathcal{H}^t(C) < \infty$. For all $K \in \N$ and $1 \leq a,b \leq  K$, we have
\begin{align*}
\mathcal{H}^t\left(\bigcup\limits_{n = 0}^{\infty}\bigcup\limits_{\mathbf{i} \in \mathcal{I}^{a + nK} }\si(C)\right) < \infty 
\iff 
\mathcal{H}^t\left(\bigcup\limits_{n = 0}^{\infty}\bigcup\limits_{\mathbf{i} \in \mathcal{I}^{b + nK}}\si(C)\right) < \infty.
\end{align*}
\end{lemma}

\begin{proof}
Fix $K \in \N$ and let
\begin{equation*}
L_{k} = \max\limits_{\textbf{i} \in \mathcal{I}^{k}} \Lip^+(\si)
\end{equation*}
for $k \in \N$. If $0 < a,b \leq K$ are distinct, we have
\begin{align*}
&\,\,\,\,\,\,\,\,\mathcal{H}^t\left(\bigcup\limits_{n = 0}^{\infty}\bigcup\limits_{\textbf{i} \in \mathcal{I}^{a + nK}} \si(C)\right) \\
&\leq \mathcal{H}^t\left(\bigcup\limits_{\textbf{u} \in \mathcal{I}^a}S_{\textbf{u}}(C) \cup \bigcup\limits_{n = 0}^{\infty}\bigcup\limits_{\textbf{i} \in \mathcal{I}^{K - b + a}}\bigcup\limits_{\textbf{j} \in \mathcal{I}^{b + nK}} S_{\textbf{i}\textbf{j}}(C)\right)\\ 
&\leq \sum\limits_{\textbf{u} \in \mathcal{I}^a}\mathcal{H}^t(S_\textbf{u}(C)) +  \sum\limits_{\textbf{i} \in \mathcal{I}^{K - b + a}}\mathcal{H}^t\left( \si\left( \bigcup\limits_{n = 0}^{\infty}\bigcup\limits_{\textbf{j} \in \mathcal{I}^{b + nK}} S_{\textbf{j}}(C)\right)\right)\\
&\leq N^aL_{a}^t\mathcal{H}^t(C) + N^{K- b + a}L_{K-b + a}^t\mathcal{H}^t\left( \bigcup\limits_{n = 0}^{\infty}\bigcup\limits_{\textbf{j} \in \mathcal{I}^{b + nK}} S_{\textbf{j}}(C)\right),
\end{align*}

by numerous applications of the scaling property of Hausdorff measure. The result follows since $a$ and $b$ were arbitrary and may be interchanged.
\end{proof}

\begin{lemma}\label{iffo}
Let $t \geq 0$ and suppose $\mathcal{H}^t(C) < \infty$. Then $\mathcal{H}^t(\mathcal{O})$ is finite if and only if 
\begin{equation*}
\mathcal{H}^t\left(\bigcup\limits_{n = 1}^{\infty}\bigcup\limits_{\mathbf{i} \in \mathcal{I}^{nK}} \si(C)\right)
\end{equation*}
is finite for some $K \in \N$.
\end{lemma}

\begin{proof}
Let $K \in \N$ and observe
\begin{align*}
\mathcal{H}^t(\mathcal{O}) &= \mathcal{H}^t\left(C \cup \bigcup\limits_{m = 1}^{K}\bigcup\limits_{n = 0}^{\infty}\bigcup\limits_{\textbf{i} \in \mathcal{I}^{nK + m}} \si(C)\right)\\
&\leq \mathcal{H}^t(C) + \sum\limits_{m = 1}^{K}\mathcal{H}^t\left(\bigcup\limits_{n = 0}^{\infty}\bigcup\limits_{\textbf{i} \in \mathcal{I}^{nK + m}} \si(C)\right),
\end{align*}
and so if
\begin{equation*}
\mathcal{H}^t\left(\bigcup\limits_{n = 0}^{\infty}\bigcup\limits_{\textbf{i} \in \mathcal{I}^{nK + K}} \si(C)\right)=\mathcal{H}^t\left(\bigcup\limits_{n = 1}^{\infty}\bigcup\limits_{\textbf{i} \in \mathcal{I}^{nK}} \si(C)\right) < \infty,
\end{equation*}
then Lemma \ref{hauskey} implies $\mathcal{H}^s(\mathcal{O})$ is also finite. The opposite implication follows by monotonicity.
\end{proof}
\subsection{Proof of Theorem \ref{hausthm}}

Let $\mathbb{I} = \{S_i\}_{i =1}^{N}$ be an IFS and $C \subseteq X$ be compact.\\
 
(i) We first note that 
$$
\mathcal{H}^t(F_\emptyset) = 0,
$$
since $s < t$ implies $\hd F_\emptyset \leq \ubd F_\emptyset \leq s < t$. Hence,
\begin{align*}
\mathcal{H}^t(F_C) &\leq \mathcal{H}^t(F_\emptyset) + \mathcal{H}^t(\mathcal{O}) 
= \mathcal{H}^t(\mathcal{O}).
\end{align*}
Thus, it suffices to show that
\begin{equation*}
\mathcal{H}^t\left(\bigcup\limits_{n = 1}^{\infty}\bigcup\limits_{\textbf{i} \in \mathcal{I}^{nK}} \si(C)\right) < \infty
\end{equation*}
for some $K \in \N$ by Lemma \ref{iffo}. Since $t > s$, it is possible to choose $K \in \N$ such that $t > s_K$ (see (\ref{skdef1}) and (\ref{skdef2})), implying
\begin{equation*}
\sum\limits_{\textbf{i} \in \mathcal{I}^{K}}\Lip^+(\si)^t <\sum\limits_{\textbf{i} \in \mathcal{I}^{K}}\Lip^+(\si)^{s_K} = 1.
\end{equation*}
It follows that
\begin{align*}
\mathcal{H}^t\left(\bigcup\limits_{n = 1}^{\infty}\bigcup\limits_{\textbf{i} \in \mathcal{I}^{nK}} \si(C)\right) &\leq \sum\limits_{n = 1}^{\infty}\sum\limits_{\textbf{i} \in \mathcal{I}^{nK}}\mathcal{H}^t\left( \si(C)\right)\\
&\leq \mathcal{H}^t(C)\sum\limits_{n = 1}^{\infty}\sum\limits_{\textbf{i} \in \mathcal{I}^{nK}}\Lip^+(\si)^t\\
&\leq \mathcal{H}^t(C)\sum\limits_{n = 1}^{\infty}\left(\sum\limits_{\textbf{i} \in \mathcal{I}^{K}}\Lip^+(\si)^t\right)^n
\end{align*}
which is a convergent geometric series and so finite, as required.\\

(ii) Suppose $\mathbb{I}$ satisfies the COSC, then
$$
\hm^t(\si(C) \cap S_\textbf{j}(C)) = 0
$$
and 
$$\mathcal{H}^t(F_\emptyset \cap \si(C))=0$$ 
for every $\textbf{i} \neq \textbf{j} \in \I^*$. Hence
\begin{align*}
\mathcal{H}^t(F_C) &= \mathcal{H}^t(F_\emptyset) + \mathcal{H}^t(\mathcal{O})\\
&= \mathcal{H}^t(F_\emptyset) + \mathcal{H}^t(C) + \sum\limits_{k = 1}^{\infty}\sum\limits_{\textbf{i} \in \mathcal{I}^k} \mathcal{H}^t(\si(C))\\
&\geq \mathcal{H}^t(F_\emptyset) + \mathcal{H}^t(C)\left(1 + \sum\limits_{k = 1}^{\infty}\sum\limits_{\textbf{i} \in \mathcal{I}^k} \Lip^-(\si)^t\right)\\
&\geq \mathcal{H}^t(F_\emptyset) + \mathcal{H}^t(C)\left(1 + \sum\limits_{k = 1}^{\infty}\left(\sum\limits_{i \in \mathcal{I}} \Lip^-(S_i)^t\right)^k\right), 
\end{align*}
and the corresponding inequality with $\Lip^+(\si)$ follows similarly. \hfill $\square$\\

\medskip
\noindent {\bf Acknowledgment}. The author was supported by a Carnegie PhD Scholarship. In addition, thanks are given to an anonymous referee for some helpful comments, and to Jonathan Fraser and Kenneth Falconer for their advice, support and discussion.

\bibliographystyle{amsplain}

\end{document}